\documentclass{amsart}
\usepackage[dvips]{graphics}
\usepackage{amsmath}
\usepackage{graphicx}

\newtheorem{theorem}{Theorem}[section]
\newtheorem{lemma}[theorem]{Lemma}
\newtheorem{corollary}[theorem]{Corollary}
\newtheorem{prop}[theorem]{Proposition}
\newtheorem{Rem}[theorem]{Remark}
\theoremstyle{definition}
\newtheorem{Def}[theorem]{Definition}

\newtheorem*{notation}{Notation}

\newcommand{\stk}[1]{\stackrel{#1}{\longrightarrow}}

\newcommand{\ra}{\rightarrow}

\newcommand{\nn}{\newline\newline}
\newcommand{\inn}{\,\in\,}

\newcommand{\N}{\mathbb N}
\newcommand{\Ss}{S_\infty}

\numberwithin{equation}{section}

\newenvironment{proofs}[1][Proof]{\textbf{#1.} }

\begin{document}

\title{On Limit Aperiodic G-Sets}

\author{Thanos Gentimis}
\address{Mathematics Department, University of Florida, Gainesville , USA }
\email{gentimis@math.ufl.edu}
\date{\today}
\begin{abstract}
We prove that the property to be limit aperiodic is preserved by the
standard construction with groups like extension, HNN extension and
free product. We also construct a non-limit aperiodic G-space.
\end{abstract}

\maketitle

\section{Introduction}
If a discrete group $G$ acts by isometries freely and cocompactly on
a metric space $X$ one can study periodic and aperiodic tilings of
$X$. A tiling of $X$ can be defined first as a tiling with one tile,
the Voronoi cell (see \cite{DS}). Using a finite set of colors one
can consider tilings of $X$ by color. Then using the "notching" one
can switch from a tiling by color to a geometric tiling. The
standard example here is $G=\mathbb Z^2$ and $X=\mathbb R^2$. Note
that the group $G$ in the above tilings is in bijection with the
tiles. Thus, construction of a geometric tiling on $X$ can be
reduced to a coloring of the group $G$. In this paper we study the
colorings of discrete groups $G$ that lead to limit aperiodic
tilings.

Let $b\in G$, a coloring $\phi$ of a group $G$ is $b$-periodic if it
is invariant under translation by $b$, i.e., for every element $g\in
G$ the elements $g$ and $bg$ have the same color. A coloring $\phi$
is aperiodic if it is not $b$-periodic for any $b\in
G\setminus\{e\}$. This can be rephrased as "The stabilizer of $\phi$
in the space of all colorings of $G$ is trivial". For infinite
groups there is a strong notion of periodicity: A coloring $\phi$ is
strongly $G$-periodic if $|Orb_G(\phi)|<\infty$. The corresponding
negation called 'weakly aperiodic' means that the orbit
$Orb_G(\phi)=G/Stab_G(\phi)$ of $\phi$ is finite. A coloring $\phi$
is called (weakly) limit aperiodic if all colorings in the closure of
the orbit $\overline{Orb_G(\phi)}$ taken in an appropriate space of
all colorings are (weakly) aperiodic.

In this paper we consider the question raised in \cite{DS}: {\em
Which groups admit limit aperiodic colorings by finitely many
colors?} This is not obvious question even for $G=\mathbb Z$. In
\cite{DS} it was answered positively for torsion free hyperbolic
groups, Coxeter groups, and groups comensurable to them.

This question can be stated in terms of Topological Dynamical
Systems theory: {\em Let $G$ be a group and $F$ be a finite set.
Does the natural action of $G$ on the Cantor set $F^G$ admit a
$G$-invariant compact subset $X\subset F^G$ such that the action of
$G$ on $X$ is free?} The dynamical system reformulation of a
corresponding question about limit weak aperiodic colorings asks
about a $G$-invariant compact subset $X\subset F^G$ such that the
orbits $Orb_G(x)$ are infinite for all $x\in X$. This was answered
affirmatively by V. Uspenskii \cite{DS}. Moreover, E. Glasner proved
that there is a minimal set $X\subset F^G$ and $x\in X$ with the
trivial stabilizer $Stab_G(x)=e$.
Despite on this progress the main question is still open. In this
paper we give a group theoretic approach. We call a group 'limit
aperiodic' (LA for short) if it admits a limit aperiodic coloring by
finitely many colors. We show that the simple group constructions
like the product, the extension, the HNN extension, and the free
product preserve the LA property. To prove these facts we introduce
the notion of LA $G$-space and prove the action theorem. In the end
of the paper we show that the main question has a negative answer
for a specific $G$-set (the natural numbers) where $G$ is a f.g. subgroup of $Aut (\mathbb Z)$.

\section{Limit Aperiodic Groups}
\begin{Def} Let $G$ be  a f.g. group. Also, let
$F$ be a finite set of elements which we can think of as colors. A
map $\phi$ from $G$ to $F$ is called a {\em coloring} of $G$.
\end{Def}
\begin{Def} Let $G$, $F$ be as above. We denote by
$F^G$ the set of all colorings from $G$ to $F$. If we consider $F$
with the discrete topology, $F^G$ with the product topology becomes
a topological space homeomorphic to the Cantor set.
\end{Def}
\begin{Def} Let $G$, $F$ as above. Then $G$ acts on $F^G$ with
the left action $\delta:G\times F^G \ra F^G$ defined by the formula
$(g\ast f)(a)=f(g^{-1}\cdot a)$ for every $g,a\inn G$ and $f \inn
F^G$.\end{Def}

Since $F^G$ is metrizable, a function $\phi$ belongs to the closure
of the orbit of $f$,  $\phi \inn \overline{Orb_G(f)}$, if and only
if $\phi=\lim\phi_k$, $\{\phi_k\}\subset Orb_G(f)$. This
is equivalent to the existence of a sequence $\{h_k\} \subset G$
with $\phi_k=h_k\ast f$. The condition $\phi=\lim(h_k\ast f)$
implies that for every $g\inn G$ there exists a $k(g)\inn \mathbb N$
with: $\phi(g)=h_k\ast f(g)$ for all $k\geq k(g)$.

\begin{Def} Let $G$,$F$ as above. A map $f:G\ra F$ is called {\em aperiodic} if
the equation $b\ast f=f$ implies $b=e$.

If the equation $b\ast f=f$ holds for some $b\in G$ we call $f$ {\em
$b$-periodic} and $b$ is called a {\em period of $f$}.
\end{Def}

\begin{Def} \textbf{(LA1)} Let $G$, $F$ be as above. A map $f: G\ra F$ is called
{\em limit aperiodic} if and only if every $\phi\inn
\overline{Orb_G(f)}$  is aperiodic.
\end{Def}

\begin{Def} \textbf{(LA2)} Let $G$, $F$ be as above. A map $f:G\ra F$  will
be called {\em limit aperiodic} if
for every $g\in G\setminus\{e\}$ there exists a finite set $S\subseteq G$,
$S=S(g)$, such that for every $h\in G$ there is a $c\in S$ with
$f(hc)\ne f(hgc)$.
\end{Def}

\begin{prop}
These two definitions are equivalent for finitely generated
groups.\end{prop}
\begin{proofs}
Suppose that $f$ satisfies the $(LA2)$ property but not the $(LA1)$.
Then there exists a $\phi\inn \overline{Orb_G(f)} $ such that $\phi$
has period $g\neq e$. Then $g^{-1}$ is also a period. Let
$\{h_k\}_{k\inn \mathbb N}\inn G$ such that $\phi = \lim h_k\ast f$.
Choose the set $S$ for that $g$. Since $S$ is finite we also have
that $g\cdot S$ is finite. From the fact that $\phi = \lim h_k\ast
f$, there exists an $n \inn \mathbb N$ such that for all $k \geq n $
and for all $x \inn S \bigcup g\cdot S$ we have $ \phi(x) = h_k \ast
f(x)$. We apply LA2 for $f$ with $g$ and $h_n^{-1}$ to obtain $c\in
S$ such that $f(h_n^{-1} c)\ne f(h_n^{-1}gc)$. This contradicts with
the fact that $g^{-1}$ is a period for $\phi$:
\[\phi(c)=(h_n\ast f)(c)=f(h_n^{-1}c)\neq f(h_n^{-1}gc)=(h_n\ast
f)(gc)=\phi(gc)=(g^{-1}\ast\phi)(c).\]

 Let's suppose now that $f$ satisfies the (LA1) but not the (LA2).
Then there exists a $g\inn G$ such that for every finite subset $S$
of $G$ there exists an $h \inn G$ with the property:
\[f(hc)=f(hgc)\] for all $c\inn S$.


Fix that $g\inn G$. Take $S_1=\{c\inn G: d(c,e)\leq 1\}$ The distance mentioned is the one induced
by the word metric in the Cayley graph of $G$.
Since $G$ is f.g. $|S_1|< \infty$, so, there exists an $h_1\inn G$ with
$f(h_1c)=f(h_1gc)$ for all $c \inn S_1$. Take
$S_2=\{c\inn G: d(c,e)\leq 2\}$. Again $|S_2|<\infty$. Then there
exists an $h_2\inn G$ such that $f(h_2c)=f(h_2gc)$ for all $c\inn
S_2$. Continue for any $k \inn \N$.

Thus we obtain a sequence $\{h_k\}_{k \inn \mathbb N}\inn G$. Taking
a subsequence we may assume that there is a limit:
\[\phi = \lim_{k\ra \infty}h_k^{-1}\ast f.\]  The claim is that
$\phi$ is periodic with period $g$. Consider an arbitrary $x\inn G$.
Name $k_1 = d(x,e)$, then  $x\inn S_k$ for all $k \geq k_1$ . Also
since $\phi$ is the limit of $h_k^{-1}\ast f$ there exists a $k_2
\inn \N$ such that for all $k\geq k_2$:
\[\phi(x)=(h_k^{-1}\ast f)(x)\] Finally since $\phi$ is the limit of
$h_k\ast f$ there exists a $k_3 \inn \N$ such that for all $k \geq
k_3$:
\[\phi(gx)=(h_k^{-1}\ast f)(gx)\]  Thus, for $k
\geq \max \{k_1,k_2,k_3\}$ we have: \[ \phi(x)=(h_k^{-1}\ast f)(x)=
f(h_kx)=f(h_k gx)=(h_k^{-1}\ast f)(gx)= \phi(gx)=(g^{-1}\ast
\phi)(x).\] Since $x$ was taken arbitrarily,  we have that $\phi$
has $g^{-1}$ as a period. This is a contradiction since $\phi$
belongs to the $\overline{Orb_G(f)}$ and $f$ has the (LA1) property.
\qed
\end{proofs}

\begin{Def} A finitely generated group $G$ will be called {\em limit aperiodic}
if it admits a limit aperiodic coloring $f:G\ra F$ with a finite set
of colors $F$.
\end{Def}
\begin{Rem} The definition of limit aperiodic groups can easily
be extended to any group and not only finitely generated ones. Both
the property (LA1) and (LA2) apply to groups without the f.g.
hypothesis. Their equivalence though depends on the fact that the
group is finitely generated. For us a group (not necessarily
finitely generated) will be limit aperiodic if it satisfies the
(LA1) property.
\end{Rem}
We recall the notion of uniform aperiodicity from \cite{DS}. Before
we introduce that notion lets establish some notation:
\begin{notation} Let $\Gamma$ be the Cayley graph of a group $G$ and $d$ be
the associated metric. We denote the {\em displacement} of $g$ at
$h$ with: \[d_g(h)=d(gh,h)\]

With $B_r(h)$ we denote the ball of radius $r$ with center $h$. Finally $\|g\|$, the norm of $g$, is the
distance between $g$ and $e$ namely: \[\|g\|=d(g,e)\]\end{notation}
\begin{Def} Let $G$ be a finitely generated group. A map $f:G\ra F$ where $F$
is a finite set (of colors) will be called {\em uniformly aperiodic}
(UA) if there exists a constant $\lambda >0$ such that for every
element $g\inn G\setminus\{e\}$ and every $h \inn G$, there exists
$b\inn B_{\lambda {d_g(h)}}(h)$ with $f(gb)\neq f(b)$.
\end{Def}
\begin{Def} A finitely generated group is called {\em uniformly
aperiodic} if there exists an $F$ and a $f$ as above, so that $f:G
\ra F$ is uniformly aperiodic.\end{Def}
$\pagebreak$

\begin{prop} If $f:G\ra F$ is uniformly aperiodic then $f$ is limit aperiodic.
\end{prop}
\begin{proofs}
We show that $f$ satisfies LA2. Let $g\inn G\setminus\{e\}$ and
$h\inn G$ arbitrary chosen. Define
 $S=B_{\lambda \|g\|}(e)$ to be the ball with center $e$ and radius
 $\lambda \|g\|$.

Clearly since $G$ is finitely generated, $S$ is finite. Assume that
there exists an $h\inn G$ such that for every $c \inn S$ we have
$f(hc)=f(hgc)$.

Denote $a=hgh^{-1}$. We apply the UA condition for $f$ with $a$ and
 $h$ to obtain $b$ in $B_{\lambda d_a(h)}(h)$ with $f(ab)\neq
f(b)$. Since $b\inn B_{\lambda {d_a(h)}}(h)$ we have that: \[d(b,h)
\leq \lambda d_a(h) = \lambda d(hgh^{-1}h,h)=\lambda d(hg,h)=
\lambda d(g,e) = \lambda \|g\| \] where the third equality comes
from the fact that the metric is left invariant.  Notice that:
\[ d(h^{-1}b,e)= d(h^{-1}b,h^{-1}h)=d(b,h)\] Thus
$d(h^{-1}b,e)\leq \lambda \|g\|$. This implies that $c=h^{-1}b$
belongs to $S$. So \[ f(b)=f(h(h^{-1}b))=
f(hc)=f(hgc)=f(hgh^{-1}b)=f(ab)\] which is clearly a contradiction.
\qed
\end{proofs}

\section{G-Sets And Limit Aperiodicity}

 The notion of limit aperiodicity can be generalized in the case of
$G$-Sets. Namely let $X$ be a space and suppose that $G$ acts on $X$
giving it the structure of a $G$-set. We will use the notation $gx=g(x)$
for $g\in G$ and $x\in X$. Fix a finite set $F$, which we can consider again as colors.

Denote by $F^X$ the set of all maps from $X$ to $F$. Then $F^X$ can
become a $G$-set under the following action: \[(g\ast f)(x)=f(g^{-1}
x)\] for all $x\inn X$, $g\inn G$ and $f \inn F^X$.

Also denote by: $Fix_G(X)=\{ g \inn G: g\cdot x = x, \forall
\,\,x\inn X\}$ the kernel of the action. $\newline$ We naturally get
the following definitions:
\begin{Def} Let $X$ be a $G$-set and let $f \inn F^X$. We call $f$
{\em limit aperiodic} if and only if for every $\phi\inn
\overline{Orb_G(f)}$ we have that $\phi$ is aperiodic meaning that
if $a\ast \phi=\phi$, then $a\inn Fix(X)$.
\end{Def}

Thus we get the definition of limit aperiodic $G$-sets:

\begin{Def} Let $X$ be a $G$-set. If there exists a finite set $F$
and a map \[f: X\ra F\] such that $f$ is limit aperiodic we say that
$X$ is a {\em limit aperiodic $G$-set}.
\end{Def}
\begin{Rem} If we consider a group $G$ acting on itself with left
multiplication then $G$ is limit aperiodic as a $G$-set if and only
if $G$ is limit aperiodic as a group, because under that action
$Fix(G)=\{e\}$.
\end{Rem}
Let $X$ be a $G$-set and let $G_x=Stab_G(x)=\{g\in G\mid gx=x\}$
denote the stabilizer of $x\in X$.
\pagebreak
\begin{theorem}\label{main} Let $X$ be limit aperiodic $G$-set and suppose
that $G$ acts transitively on $X$. Fix $x\inn X$ such that
$X=Orb_G(x)$. If $Stab_G(x)$ is a limit aperiodic group for some
$x\in X$ then $G$ is a limit aperiodic group. \end{theorem}

\begin{proofs}
Let $\phi:X\ra F_1$ be a limit aperiodic map for $X$ and let
$\psi:G_x\ra F_2$ be a limit aperiodic map for $G_x$. We know that
there exists a bijection $\pi$ between the set of left cosets
$G/G_x$ and the orbit $Orb_G(x)$. Fix a set of representatives in
$G$ namely $\{a_i: i \inn I\}$ for the quotients $G/G_x$. Then
$\pi(a_jG_x)=a_jx$. $\nn$ Define $f=(f_1,f_2):G\ra F_1\times F_2$ by
$f_1(g)=\phi(gx)$ and $f_2(g)=\psi(a_j^{-1}g)$ where $g\inn a_jG_x$.
We will prove that this $f$ is a limit aperiodic map. Suppose that
this is false. Then there exists a map
$\overline{f}\in\overline{Orb_G(f)}$ and an element $a \inn G$ such
that $\overline{f}$ has $a$ as a period, i.e., $(a\ast
\overline{f})=\overline{f}$ for all $x\inn X$. Let
$\overline{f}=\lim h_k \ast f$ where $h_k \inn G$

$\nn$ \underline{Case 1)} Suppose $a\,\notin \,G_x$. Consider the
limit
\[\overline{\phi} = \lim_k h_k \ast \phi.\] Note that we can always choose a
subsequence of $h_k$ such that the limit exists. For convenience we
keep the same indices for the subsequence. Clearly,
$\overline{\phi}\in\overline{Orb_G}(\phi)$. Given $g\inn G$, there
exists  $k_0\inn \mathbb N$ such that for all $k\geq k_0$ we have:
\[\overline{\phi}(gx)=(h_k\ast\phi)(gx)\] Also there exists a
$k_1\inn \mathbb N$ such that for all $k\geq k_1$ we get:
$\overline{f}(g)=(h_k\ast f)(g)$. Let $k_2=\max\{k_0,k_1\}$ then for
all $k\geq k_2$ we have $\overline{f}(g)= f(h_k^{-1}g)$. Therefore,
\[\overline{f_1}(g) =
f_1(h_k^{-1}g)=\phi(h_k^{-1}gx)=(h_k\ast\phi)(gx)=\overline{\phi}(gx).\]
Following the previous proof and replacing $g$ with $a^{-1}g$ we
find a $k_3\inn \mathbb N$ such that for all $k \geq k_3$ we have:
\[\overline{f_1}(a^{-1}g)=\overline{\phi}(a^{-1}gx)\]
Since $\overline{f}_1$ has period $a$ we have $(a\ast
\overline{f})(g)=\overline{f}(g)$. Hence, \[(a\ast
\overline{\phi})(gx)=\overline{\phi}(a^{-1}gx)=
\overline{f_1}(a^{-1}g)=(a\ast
\overline{f_1})(g)=\overline{f_1}(g)=\overline{\phi}(a^{-1}gx).\]

Since $g$ is arbitrary and $G$ acts on $X$ transitively we get that
for every $y \inn X$, $(a\ast \overline{\phi})(y) =
\overline{\phi}(y)$. Since $\phi$ is limit aperiodic we have that $a
\inn Fix(X)$. But: \[Fix(X) = \bigcap_{s\inn S} Stab_G(s)\] Thus, $a
\inn Fix(X)\subseteq Stab_G(x)=G_x$ contradiction.
$\nn$\underline{Case 2)} Suppose that $a\inn G_x$. Let
$\{h_k\}$ be a sequence of elements of $G$ such that $h_k^{-1}$
belongs to the coset $a_kG_x$. Thus $\delta_k=h_ka_k$ belongs to
$G_x$.  Taking a subsequence we may assume that there are the limits
\[ \overline {\psi} = \lim_k \delta_k\ast\psi \ \ \ \ \ \text{and}\ \ \ \
\overline {f} = \lim_k h_k\ast f.\]  Notice that $\overline {\psi}
\inn \overline {Orb_{G_x}(\psi)}$. Let $h\inn G_x$. Then there
exists a $k_0\inn \mathbb N$ such that for all $k\geq k_0$ we have
\[\overline{\psi}(h)=(\delta_k\ast\psi)(h)\]
Also there exists a $k_1\inn \mathbb N$ such that for all $k \geq k_1$ we get:
\[\overline{f}|_{G_x}(h)=(h_k\ast f)(h)\]
Then for all $k \geq \max\{k_0,k_1\}$ we have:
$\overline{f}|_{G_x}(h)= f(h_k^{-1}h)$.
Hence, \[(\overline{f}|_{G_x})_2(h) = \psi(a_k^{-1}h_k^{-1}h)
=\psi((h_ka_k)^{-1}h)=((h_k a_k)\ast \psi)(h)=(\delta_k\ast
\psi)(h)=\overline{\psi}(h).
\]
 Notice now that:
\[a\ast
\overline{\psi}=a\ast\overline{f}|_{G_x}=\overline{f}|_{G_x}=\overline{\psi}.\]
Thus $\overline{\psi}$ is periodic. Contradiction since $\psi$ is
limit aperiodic. This concludes the proof.\qed
\end{proofs}

The following is obvious.
\begin{lemma}\label{obvious} If $X$ is a $G$-space and $H$ is any group
acting on $X$ such that the action of $H$ factors through
the action of $G$. Then if $X$ is a limit aperiodic $G$-space then
it is also a limit aperiodic $H$-space.
\end{lemma}

\begin{corollary} If $G$ and $H$ are limit aperiodic groups and
\[1\ra G \stk{\tau} E \stk{\pi} H \ra 1\] is a short exact sequence then $E$ is
also limit aperiodic. \end{corollary}
\begin{proofs}
Obviously $E$ acts transitively on $E/G=H$ with left multiplication.
By Lemma \ref{obvious} $H$ is a limit aperiodic $E$-space. Note that
$Stab_E(e_H)=G$ is limit aperiodic. If we apply Theorem \ref{main}
we get the corollary.\qed\end{proofs}

Obviously we obtain the following
\begin{corollary} If $G$ is limit aperiodic and $H$ is limit
aperiodic then $G\times H$ is limit aperiodic.
\end{corollary}
\begin{corollary} If $H$ is a limit aperiodic group and $\theta:H\ra
H$ is a group automorphism then the HNN extension $\star_\theta H$
is limit aperiodic.
\end{corollary}
\begin{proofs}  We know that if $H=<S|T>$ where $S$ is a set of
generators and $T$ is a set of relations then $G=\star_\theta H$
admits the following presentation \[\star_\theta H=<S,t|T\cup
\{t^{-1}xt=\theta(x), x\in S\}>.\]  Note that $G$ acts transitively
on the set $G/H$ of all left cosets of $H$. Note that
$G/H=\{t^iH\mid i\in\mathbb Z\}\cong\mathbb Z$. Thus, $G$ acts on
$\mathbb Z$ by translations with $Stab_G(\{H\})=H$ \cite{Serre}. We
color the $\mathbb Z$ with a Morse-Thue sequence as in \cite{DS},
i.e. $\phi:X\ra \{0,1\}$ with $\phi(t^iH)=m(|i|)$ where $m:\mathbb N
\ra \{0,1\}$ is the Morse-Thue sequence \cite{Morse},\cite{Thue}. As
it shown in \cite{DS} this map is limit aperiodic with respect. By
Lemma \ref{obvious} it is a limit aperiodic $G$-set. Theorem
\ref{main} completes the proof. \qed
\end{proofs}

In order to prove the fact that the free product of limit aperiodic
groups is limit aperiodic we need the following notions. Let $A$ and
$B$ be two groups. We will construct a set $X$ such that the free
product  $G=A\star B$ acts on $X$ freely and transitively and $X$ is
a limit aperiodic $G$-space. Let $T_0$ be the Bass-Serre tree
associated with $A\star B$. We recall that the vertices of $T_0$ are
left cosets $G/A\cup G/B$ and the vertices of the type $xA$, $xB$
and only them form an edge $[xA,xB]$ in $T_0$. Thus the edges of
$T_0$ are in the bijection with $G$. Note that $G$ acts on $T_0$ by
left multiplication. Let $T$ be the barycentric subdivision of
$T_0$ and let $X$ be the set of the barycenters (of edges). We will
identify the tree $T$ with the set of its vertices. Then the group
$G$ acts by isometries on $T$ yielding three orbits on the vertices
$X=Orb_G(e)$, $G/A$ and $G/B$. We regard $T$ as a rooted tree with
the root $e$. Let $\|x\|=d_T(x,e)$ denote the distance to the root.

\begin{lemma}\label{T} Let A,B be limit aperiodic groups and let $G=A\star B$,
$X$, $T$ defined as above. Then $T$ is a limit aperiodic $G$-set.
\end{lemma}
\begin{proofs}
Let $\pi: G\ra A$ with $\pi
(w)=\pi(a_1b_1a_2b_2...a_nb_n)=a_1a_2...a_n$ and $\theta: G\ra B$
with $\theta(w)=\theta(a_1b_1a_2b_2...a_nb_n)=b_1b_2...b_n$. Clearly
both $\pi$ and $\theta$ are group homomorphisms. Let also $f_A:A\ra
F_A$ be the limit aperiodic map for the group $A$ and $f_B:B\ra F_B$
be the limit aperiodic map for group $B$. Also let $\nu:Z\ra
\{0,1,2\}$ be the variation of the Morse-Thue sequence which has no
words $WW$ (see for example \cite {DS}). Also fix $e$ to be the
vertex representing the identity element in $T$. Then consider a coloring of $T$ as
follows:
\[f:T\ra \{0,1,2\}\times \{0,1,2\} \times (F_A\bigcup\{\alpha\})\times
(F_B\bigcup\{\beta\})=F\] where $f:=(f_0,f_1,f_2,f_3)$ with:
$f_0(x)=\nu(\|x\|)$, $f_1(x)=\|x\|mod3$, $f_2(x)=f_A(\pi(x))$ if
$x\inn X$ and $f_2(x)=\alpha$ if $x\inn T-X$. Finally let
$f_3(x)=f_B(\theta(x))$ if $x\inn X$ and $f_3(x)=\beta$ if $x\inn
T-X$. The group $G$ acts on the space of colorings $F^T$ as follows:
\begin{eqnarray*} (g\ast f)(x)&=&((g\ast f_0)(x),(g\ast f_1)(x),
(\pi(g)\ast f_2)(x),(\theta(g)\ast f_3)(x)\\ &=&
(f_0(g^{-1}x),f_1(g^{-1}x),f_2(\pi^{-1}(g)x),f_3(\theta^{-1}(g)x)).
\end{eqnarray*}
Suppose that $f$ is not a limit aperiodic map. Then there exists a
coloring \[\psi=(\psi_0,\psi_1,\psi_2,\psi_3)\] such that $\psi\inn
\overline{Orb_G(f)}$ and $\psi$ has a period $b\inn G\setminus
Fix(T)$. Let $\psi=\lim g_k\ast f$. Then $\psi_A=\lim\pi(g_k)\ast
f_A$ has period $\pi(b)$. Indeed, for every $x\in A\subset G\cong X$
for large enough $k$,\begin{eqnarray*} (\pi(g_k)\ast
f_A)(x)=f_A(\pi(g_k^{-1}x))=f_2(g_k^{-1}x)=(g_k\ast f_2)(x)=(g_k\ast
f_2)(bx)\\
= f_2(g_k^{-1}bx)=f_A(\pi(g_k^{-1}bx))=(\pi(g_k)\ast
f_A)(\pi(b)x).\end{eqnarray*} Similarly $\psi_B=\lim \theta(g_k)\ast
f_B$ has a period $\theta(b)$. Thus, $\pi(b)=e_A$ and
$\theta(b)=e_B$.

Notice that $(\psi_0,\psi_1)\inn \overline{Orb_G(f_0,f_1)}$. Denote
$\xi=(\psi_0,\psi_1)$ and $\phi = (f_0,f_1)$. Then $\xi$ is a
coloring of a simplicial tree $(T)$ on which $G$ acts by isometries.
Moreover $\xi\inn \overline{Orb_G(\phi)}$. From proposition $4$,
page 318 in \cite{DS} we have that $b\ast\xi\neq \xi$ for all $b\inn
G$ with unbounded orbit $\{b^kx_0|k\inn \mathbb Z\}$. This clearly
implies that $b\ast\psi \neq \psi$ for every $b\inn G$ with
unbounded orbit. On the other hand $\psi$ has period $b$ and thus we
have $\{b^kx_0|k\inn\mathbb N\}$ is bounded. This implies that $b$
fixes a point in $T$. Call that point $x_1$. Since the action of $G$
on $X$ is free, $x_1\notin X$. Thus, $x_1\in G/A$ or $x_1\in G/B$.
Assume the later, $x_1=wB$ for some $win G$. Since $b$ fixes $wB$,
$b=wb'w^{-1}$ for some $b'\in B\setminus\{e_B\}$. Then
$\theta(b)=\theta(w)b'\theta(w)^{-1}\ne e_B$. Contradiction. \qed
\end{proofs}

\begin{lemma}\label{X} Let $G=A\star B$, $T$, $X$, $f$, $F$ as above. Then
$X$ is a limit aperiodic $G$-set.
\end{lemma}
\begin{proofs} Note that in the rooted tree $T$ every vertex $x\ne e$ has
a unique predecessor denoted $pred(x)$. We define $f':X\ra F\times
F$ as $f'(x)=(f|_X,\hat f)$ where $\hat f(x)=f(pred(x))$. We show
that $f'$ is limit aperiodic.

Suppose that $f'$ is not limit aperiodic. Then there exists a
sequence $\{g_k\}\inn G$ s.t.
\[\psi'=\lim_kg_k\ast f'\] and $b\inn Fix_G(X)$ with $b\ast
\psi'=\psi'$. We may assume that there is the limit $\psi=\lim
g_k\ast f$. In view of Lemma \ref{T} it suffices to show that $\psi$
is $b$-periodic.
It is $b$-periodic on $X$, so it suffice to check that it is
$b$-periodic on $T\setminus X$. Let $z\in T\setminus X$. We check
that $\psi(bz)=\psi(z)$. Since the root $e$ lies in $X$, we may
assume that $z\ne e$. Let $x_0=pred(z)$ and let $x_1$ be such that
$z=pred(x_1)$. We note that $x_1$ is not unique. So we fix one. Note
that $x_0,x_1\in X$. There is $k_0$ such that for $k\ge k_0$
\[ \psi'(x_i)=(g_k\ast f')(x_i),\ \ \ \ \psi'(bx_i)=(g_k\ast
f')(bx_i),\ \ \ i=0,1\] and

\[ \psi(z)=(g_k\ast f)(z),\ \ \ \ \psi(bz)=(g_k\ast
f)(bz).\] Fix $k\ge k_0$. Since $G$ acts on $T$ by isometries the
distance from $g_k^{-1}z$ to $g_k^{-1}x_i$, $i=0,1$ equals 1. There
are three possibilities:

\[(1)\ \ \ \ \ \ \  g_k^{-1}x_0<g_k^{-1}z<g_k^{-1}x_1,\]

\[(2)\ \ \ \ \ \ \ g_k^{-1}x_1<g_k^{-1}z<g_k^{-1}x_0,\] and

\[(3)\ \ \ \ \ \ \ g_k^{-1}z<g_k^{-1}x_i,\ \ \ i=0,1.\]

We apply $g_k^{-1}bg_k$. In view of the fact that
$f_1(g^{-1}_kx_i)=f_1(g_k^{-1}bx_i)$, $i=0,1$ we obtain

\[g_k^{-1}bx_0<g_k^{-1}bz<g_k^{-1}bx_1\] in the case (1) and

\[g_k^{-1}bx_1<g_k^{-1}bz<g_k^{-1}bx_0,\] in the case (2).
Then in the case (1)
\[ f'(g_k^{-1}bx_1)=(g_k\ast
f')(bx_1)=\psi'(bx_1)=\psi'(x_1)=(g_k\ast
f')(x_1)=f'(g_k^{-1}x_1).\] Hence
$f_0(g_k^{-1}bx_1)=f_0(g_k^{-1}x_1)$. Therefore
\[f(g_k^{-1}bz)=f(g_k^{-1}z).\] Thus,
\[\psi(bz)=(g_k\ast f)(bz)=f(g_k^{-1}bz)=f(g_k^{-1}z)=(g_k\ast
f)(z)=\psi(z).\] In the case (2) we consider $x_0$ instead of $x_1$.

In the case (3) $g_k^{-1}bz$ is the predecessor of either
$g_k^{-1}bx_0$ or $g_k^{-1}bx_1$ (or both). Assume the first. Then
from the $b$-periodicity of $\psi'$ it follows that
$f_0(g_k^{-1}bx_0)=f_0(g_k^{-1}x_0)$. Since
$f_0(g_k^{-1}bx_0)=f(pred(g_k^{-1}bx_0))$ and
$f_0(g_k^{-1}x_0)=f(pred(g_k^{-1}x_0))$, we obtain
$\psi(bz)=f(g_k^{-1}bz)=f(g_k^{-1}z)=\psi(z)$. \qed
\end{proofs}
\pagebreak

\begin{theorem}\label{amalgam} Let $A$,$B$ be limit aperiodic groups.
Then $G=A\star B$ is
a limit aperiodic group.
\end{theorem}
\begin{proofs} By Lemma \ref{X} $X$ is a limit aperiodic $G$-set.
Note that $G$ acts on $X$ as above transitively, and $Stab_G(x_0)=\{e\}$
is a limit aperiodic group. By Theorem \ref{main} we have that $G$
is a limit aperiodic group.\qed
\end{proofs}
We finish this paper with an example of a $G$-set which is not limit aperiodic.

\begin{prop} Consider the automorphism group of the integers $Aut(\mathbb Z)$. Let $s:\mathbb Z\ra \mathbb Z$ with $s(n)=n+1$ and
$t:\mathbb Z\ra \mathbb Z$ with $t(0)=1$, $t(1)=0$ and $t(n)=n$ if $n\neq 1$ and $n\neq 0$. Let $S=<s,t>$ then $\mathbb N$ is not a limit aperiodic $S$-set.
\end{prop}
\begin{proofs} Suppose that $\N$ was a limit aperiodic
$S$-set. Then let $F$ be a set of colors with $|F|<\infty$ and a
map $f:\mathbb N \ra F$ such that $f$ is a limit aperiodic map under
the action of $S$. Since $|F|<\infty$ and $|\N|=\infty$ there
exists at least one $a\inn F$ such that infinitely many $a_n\inn \N$ have $f(a_n)=a$. Choose a strictly
increasing sequence in $\N$ such
that $f(a_n)=a$ for all $n\inn \N$. Consider the
following elements in $S$: \begin{eqnarray*} h_1 &=& (1,a_1) \\
h_2 &=& (1,a_1)(2,a_2) \\ \dots \\ h_n&=&(1,a_1)(2,a_2)\dots
(n,a_n)\\ \dots
\end{eqnarray*} where $(i,a_i)$ is the transposition that
takes $i$ to $a_i$. With $s$ and $t$ we can construct all the transpositions. Thus all $a_i$ belong to $S$.
Clearly if $n\geq k$, $n,k \inn \N$ we have that
$h_n(k)=a_k$. Consider now the sequence $\{h_n^{-1}\ast f\}$ and take
a converging subsequence. For convenience in notation let us keep
the same indices for the subsequence. Thus if:
\[\psi=\lim_{n\ra \infty} h_n^{-1}\ast f\] we have that $\psi\inn
\overline{Orb_{\Ss}(f)}$. Thus $\psi$ has to be aperiodic. The claim
is that $\psi(k)=a$ for all $k\inn \N$ and thus $\psi$ is clearly
periodic which will lead to a contradiction. This is easy to see
since let $k\inn \N$. Then for that $k$ there exists an $n_1$ such
that for all $n\geq n_1$ we have that $\psi(k)=(h_n^{-1}\ast f)(k)$.
Thus for $n= \max\{k,n_1\}$ we have that:
\[\psi(k)=(h_n^{-1}\ast f)(k)=f(h_nk))=f(a_k)=a.\] This concludes the
proof.\footnote{ When this paper was finalized we received information about a preprint on limit aperiodic colorings of groups by Gao, Jackson and Seaward
found in http://www.cas.unt.edu/~sgao/pub/pub.html} 
\end{proofs}


\begin{thebibliography}{99}
\bibitem{DS}
A.Dranishnikov and V.Schroeder {\em Aperiodic Colorings And Tilings
Of Coxeter Groups} Groups Geom. Dyn. 1(2007), 311-328.

\bibitem{Serre}
J.P. Serre {\em Trees}, Springer-Verlag(2003), ISBN 3540442375

\bibitem{Morse}
H.M. Morse and G.A. Hedlund {\em Unending chess, symbolic dynamics
and a problem in semigroups} Duke Math. J. 11 (1944), 1-7.

\bibitem{Thue}
A. Thue {\em Uber unendliche Zeichenreihen} Christiania Vidensk.
Selsk. Skr. No 7 (1906),  p. 22.

\end{thebibliography}
\end{document}